\documentclass{amsart}
\usepackage{amsmath}
\usepackage[pdftex,colorlinks]{hyperref}

\newcommand{\dt}{\mathrm{d}t}

\begin{document}
\title{Fractions in Elementary Education}
\date{Version 1.2, November 2013}
\author{Frank Quinn}
\maketitle
\section*{Introduction}   
This paper is one of a series in which elementary-education  practice is analyzed by comparison with the history of mathematics,  mathematical structure, modern practice, and (occasionally) cognitive neuroscience. The primary concerns are: Why do so many children find elementary mathematics difficult? And, why are the ones who succeed still so poorly prepared for college material needed for technical careers? 
The answer provided by conventional wisdom is essentially that mathematics is difficult. Third-graders are not developmentally ready for the subtlety of fractions, for instance, and even high-performing students cannot be expected to develop the skills of experienced users. However we will see that this is far from the whole story and is probably wrong: elementary-education fractions are genuinely harder and less effective than the version employed by experienced users. Experts  discard at least 90\% of what is taught in schools. Our educational system is actually counterproductive for skill development, and the reasons for this are an important secondary concern.

\subsection*{History} In a nutshell, what we think of as mathematics originated in Greece about 2400 years ago and developed slowly over the following 2000 years. In the 1600s, driven by the mathematization of science, professional development changed directions and accelerated considerably. This is when ratios were replaced by fractions, for instance, and when continuous magnitudes were finally decomposed into real numbers and units. Elementary education is still heavily infuenced by the methodology of the 1500s. In other words, what we find is that children are taught methods and perspectives that professionals  have considered obsolete for 400 years. It would be strange if children were \emph{not} confused, and if their early training was \emph{not} a barrier to learning modern material. 

As for why this disconnect persists,  educators reject skill development as a primary goal on the grounds that most people have little use for mathematical skills later in life. There is some sense to this, but children interested in technical professions \emph{will} need these skills and, for better or worse, responsibility for getting this done is located in the schools. Unfortunately the nature of the subject has enabled a failure of accountability. If this were music we could say that picking out melodies on a piano is ok for ``understanding'', but students interested in careers in music are not well served by four years of it. By performance standards they would sound  bad, and it would be obvious that they would have to unlearn a lot of bad habits before they could develop proficiency. But it takes a lot of training to ``hear'' off-key mathematics, so instead of actual performance many educators now focus on \emph{philosophy} of performance. It must be said that their philosophy is very attractive. And, since western philosophy is strongly influenced by  ancient mathematical methodology, the ancient methodologies fit right in.  As for why it persists,  people at the college level who are supposed to teach melody-pickers how to play science and engineering do complain, but these complaints are incoherent and easily deflected by the attractive philosophy. 
\subsection*{The Common Core Standards}
The detailed analysis focuses on the treatment of fractions in the  US Common Core State Standards in Mathematics \cite{ccssm}. A micro-level analysis requires a specific context, and this  seems to be the best choice.  For one thing, extensive debates and revisions in the development process ensure that they approximate a consensus view.  In contrast, the published literature is essentially a repository of opinions  without coherence or closure mechanisms. Stability is another virtue: the  adoption process was difficult enough to make change in the near future unlikely. In contrast, individual state standards, texts, etc.~sometimes change substantially on short notice. Finally, the linkage between literature, standards, texts, tests, etc.~is rather weak. For example the Common Core standards has a ``Sample of Works Consulted''  but it is not used in a scholarly way: there are no citations in the body of the document\footnote{The Glossary is somewhat better about this.}. Are these suggestions for further reading? Were some followed and others rejected? In any case analyzing references would not illuminate the document. 
\subsection*{Outline} in the first section five problematic features of fractions in the Common Core are identified and analyzed. The second section  describes modern fractions. The standard question ``what are fractions?'' evokes   technical and unhelpful answers so
instead we observe that  after the 1600s education stayed largely constant while expert practice changed. Identifying the driver of change in expert practice enables us to identify what is ``better'' about the modern version, and give a description that puts this special virtue up front. This is where the basic simplicity emerges. 

The final section addresses the secondary question ``why didn't elementary education eventually follow the experts?''

 \section{Elementary-education fractions}
We analyze the treatment of fractions in the Common Core from their introduction in grade three through their use in proportion problems in grade eight. Five problematic features are identified. Four come directly from the text: parts-of-a-whole,  visual fraction models,  word problems, and  ratios  and proportions. The fifth feature is the overall incoherence and dysfunctionality of the material.

\subsection{Parts of a whole}\label{ssect:ptsofwhole} The dominant elementary-education view is that `fraction' is an amalgam of sub-constructs\footnote{In both language and perspectiven the K-8 Common Core standards follow the analysis of Kieren \cite{kieren}, with later refinements (cf.~\cite{char07}). The subconstructs are `parts-of-a-whole', `ratio', `operator', `quotient', and `measure',
though some of the interpretations of `measure'  are outside the K-8 window.}.
 `Parts-of-a-whole' is the most important of these and, as befits a primary concept, the CCSS-M grade 3 section on  Number and Operations  begins with it:\medskip
\begin{quote}1. Understand a fraction \( 1/b\) as the quantity formed by 1 part when a
whole is partitioned into \(b\) equal parts\dots
\end{quote}
Further down we find:
\begin{quote}3d. Compare two fractions with the same numerator or the same
denominator by reasoning about their size. Recognize that
comparisons are valid only when the two fractions refer to the
same whole.\end{quote}
There are two further references to `parts of a whole' in grade 3, three in grade 4 and then, after one last caution in grade 5, the term disappears. What is the function of this term, and why does it become unnecessary in later grades? 
\subsubsection{The precision-avoiding loophole}\label{sssect:loophole} 
In the preamble to the grade 3 standards we find:
\begin{quote}
Students understand that the size of a
fractional part is relative to the size of the whole. For example, 1/2 of the
paint in a small bucket could be less paint than 1/3 of the paint in a larger
bucket\dots\end{quote}
But 1/2 is greater than 1/3 no matter how big the buckets are.   This is a difficulty with physical applications, not fractions, and 
``different wholes'' is used to avoid having to be clear about this distinction. 

The amount of paint in a bucket is a bit abstract. Grade 3 examples typically refer to visible comparisons that connect with  `visual fraction models',  so we use one of these for detailed analysis.
 \begin{quote}Al's pizza was divided into 8 pieces and he got 2. Jose's pizza was divided into 6 pieces and he got 1. Who got the most pizza?\end{quote}
 In the context of eating flat things, ``most'' is taken to refer to area. This is consistent with naive visual comparisons of plane figures by area. The data given, however, concerns vertex angle rather than area, and if the pizzas have different diameters then area won't correlate with vertex angle.   The ancient approach is to declare that analytic comparisons in such cases don't make sense because they refer to different ``wholes''. This avoids the  issue because for a single ``whole'', vertex angle will correlate with most \emph{any} measure: amount of pepperoni, thickness of crust, etc., as well as the visual area. Pedagogically convenient but unsatisfactory for later work. 

The modern perspective is that the formulation of the problem is sloppy. The data concerns vertex angle but the  answer implicitly concerns area. One way to make it honest is to specify the criterion for the answer: ``who got the most \emph{vertex angle} of pizza?''  This sounds silly, not least because there is a disconnect between the explicit answer criterion and our implicit interest in area. The version ``who got the most \emph{area} of pizza'' addresses our implicit interests but exposes the disconnect between data and implicit question. 
Perhaps this is just a bad problem. In any case it should be clear that data/question disconnect should be addressed with more precision or more data, and dodging it with a ``different wholes'' restriction on the fraction concept is inappropriate. 

 Other uses of the `different wholes' loophole are described in the section on ratios and proportions. The parts-of-a-whole picture fades away in grade 5 because students are using a more capable version of ratios and the loophole is no longer needed, see \S\ref{ssect:oldratios}. 

\subsubsection{Ratios, not fractions} `Parts-of-a-whole' describes a primitive form of ratios that differs from fractions in several ways. The first difference is the 2400-year-old perspective that a ratio  \(a\):\(b\) is a \emph{pair} of things (here `part' and `whole'), and these exemplify a relationship.  Ratio calculations proceed by manipulating pairs in ways that  preserve the relationship, see \S\ref{ssect:oldratios} for details.  The more recent fraction perspective is that \(\frac{b}{a}\) denotes a single rational thing  that encodes the relationship in a particular way, and calculations use fraction arithmetic. 

The other difference that is significant here is that when the definition of ratios is made mathematically precise, they turn out to correspond to angles in polar coordinates, see  \S\ref{ssect:carefulratios}. `Parts-of-a-whole'  makes perfect sense for angles: the circle is the whole, and parts correspond to segments. The phrase honestly reflects the mathematical structure of ratios, but shows they are \emph{different} from fractions rather than the same. This difference is clearly visible in pie-chart visual fraction models, see \S\ref{sssect:piecharts}.

   \subsubsection{Proper fractions} 
The phrase `parts-of-a-whole' is often understood as implying that fractions should not be greater than 1. This reading is legitimate because it is essentially correct for ratios. It is also embedded in terminology such as `proper' and `improper' fractions, and the `mixed-number' practice of expressing large fractions as integers plus proper fractions. All this makes sense for ratios but is not helpful for fractions.   When educators try to explain how a part can be bigger than the whole, for instance, they are not only fudging the literal meaning and logical consequences of `parts-of-a-whole', but doing something that their own terminology identifies as ``improper''. 

 \subsubsection{Other sources} 
Historically the phrase `parts-of-a-whole' may also have specified a particular form of ratio. For instance, in ``mortar mix consists of three parts sand to one part cement'', the numbers 3 and 1 refer to parts of a \emph{sum}, not of a whole.  The corresponding parts-of-a-whole form would be ``mortar mix consists of sand, three parts in four; and cement, one part in four.''  Note that this ambiguity does not occur in rate and fraction formulations. 
 \subsubsection{Summary} `Parts-of-a-whole' 
is a legacy from early treatments of ratios, and is used primarily to dodge difficulties with over-simplified word problems. The ratio aspect is ironic because the heading of the CCSS-M section where it first occurs is ``Develop understanding of fractions as numbers''. It also illustrates the very low standards of precision in education. But is it really a problem, or just a  pedagogically convenient fudge that will be corrected  before anyone notices? Unfortunately, in this case the logical inconsistencies are close to the surface, and in fact the `parts-of-a-whole' view of ratios was obsolete long before ratios themselves were phased out. Students can be taught that they are the same anyway, but this develops confusions and mindsets that are hard to overcome. Remember that it took professionals 2000 years to go from ratios to modern fractions, and that the ratio mindset blocked the use of negative numbers in the West for centuries.

\subsection{Visual fraction models} \label{ssect:visualmodels}
In section 3NF of the Common Core we find:
\begin{quote}3b. Recognize and generate simple equivalent fractions, e.g., 1/2 =
2/4, 4/6 = 2/3. Explain why the fractions are equivalent, e.g., by
using a visual fraction model.\end{quote}
The phrase ``visual fraction model'' occurs twice more in grade 3, six times in grade 4, eight times in grade 5, twice in grade 6, and then vanishes. There are also references to the general practice that do not use this exact phrase. 

The phrase is usually connected with ``show'', ``explain'', or ``understand'', and in 11 of the 17 instances this is the only approach mentioned.  Visual models are plainly the preferred approach, at least through grade 5.  What does ``visual fraction model'' mean, what are the consequences of its use, where does it come from, and why is it not used in higher grades? 

\subsubsection{Mathematics based on vision} Precise logical reasoning is the core activity in mathematics, but reasoning has to be \emph{about} something. The practice developed by the Greeks 2400 years ago is to describe objects by example, and infer properties from physical experience and perception. We are then supposed to reason logically with these inferred properties.

``Visual fraction models'' is an instance of this ancient Greek practice. In this case  visual impressions of area are used to infer what happens when you cut a region into pieces. The Common Core also has a few references to use of `number line' models. In principle these would use visual comparisons of lengths of line segments. Human perceptions of length are weaker than perceptions of  area, however, so it is common to thicken the line into a ribbon and compare areas of pieces of the ribbon. 

Six of the 17  Common Core references to visual fraction models include ``and/or equations''. This addendum clarifies that though visual models are plainly preferred, more-symbolic approaches are not ruled out.   The equation approach was well-established in the Arabic world by 900, and in Europe in the 1200s. Methods needed for fully symbolic modern treatments were available by 1600. In other words visual fraction models may be useful for illustration, but as a main approach they have been obsolete for 400 to 1000 years. In detail,
\begin{enumerate}\item Some aspects of the visual approach are mathematically wrong (see below), so the presentation has to be vague and ambiguous.
\item it is  too clumsy to help with computation. 
\item Accuracy of human vision limits the method to small integers, and in the Common Core this shows up as severe limits on the size of denominators allowed in grade 4, for instance.
\item It is a very poor model for mathematical reasoning. 
\end{enumerate}
  Modern practice is to  define objects precisely, and in ways that allow properties to be extracted by logical reasoning. The drawback is that some logical work is necessary before anything can be done with such objects, whereas the visual approach has almost no learning curve. In other words, the visual product is bad, but educators like it because it is cheap.
\subsubsection{Angles and pie charts}\label{sssect:piecharts} Pie charts are the most common visual fraction model.  The information presented in a pie chart concerns vertex angles (segments of a circle); enlarging these into pie slices makes the information accessible to visual perceptions of area. We saw in \S\ref{sssect:loophole} that the difference between angle and area is itself a source of confusion, but we have other concerns here. 

On a heuristic level pie charts connect very well with  `parts-of-a-whole'. The circle is the whole, and parts correspond to intersections of pie slices with the circle (vertex angles). This also fits well with the mixed-number form for improper fractions. Think of the circle as the real numbers modulo 1. Any real number can be written as an integer plus a  ``proper''  part \(t\) with \(0\leq t<1\).  The proper  part is a representative of the number modulo 1. Telling time with an analog clock uses the same pattern: the clock records elapsed time modulo 12 hours, or in other words as proper fractions of a half-day.  
The problem is that angles correspond to ratios---more precisely \emph{oriented} ratios, see \S\ref{ssect:carefulratios}---not fractions. 

Subdividing a rectangle into pieces is a visual model that does represent fractions. For one thing it inverts the stacked-rectangle model for multiplication. Also, rectangles are not seen as representing a unit area (squares do that), so there is no parts-of-a-whole aspect and improper fractions can be sensibly represented. 

The differences in the structures underlying  rectangular and pie-chart models are pretty clear. 
If students are confused by explanations of why they are somehow the same, maybe the reason is that they are not the same. 

\subsection{Word problems, ratios, and proportions} \label{ssect:oldratios}  
In the K--8 part of the Common Core document, the phrases ``word problem'' or  ``real-world'' occur a total of 51 times.  In the fraction material there are 22 sample problems (set in italics), and all but one of these concern physical situations. Illustrations of arithmetic operations in Tables 1 and 2 in the Glossary are all word problems.  Word problems have always had a central role in mathematics education but the emphasis in the Common Core  is particularly heavy. This emphasis has far-reaching consequences but here we restrict to consequences for the study of fractions. One of these is that it sustains the use of ratios. 

Ratios and proportions are obsolete precursors of fractions whose main virtue is that they give a primitive way to handle some of the physical aspects of word problems. They are
 major topics in grades 6 and 7 in the Common Core. They are not mentioned before grade 6 even though implicit  use  begins in grade 3, and they essentially disappear in grade 8. We will see that fraction arithmetic has replaced the original ratio manipulations as the ``computational engine'', so in practice ratios serve only as a conceptual overlay to deal with word problems. In grade 8 students begin to use the  modern (ca.~1650) calculus of units, which gives a much more effective conceptual overlay. 
 
 \subsubsection{Sharing rice}
To illustrate the issues we work through an example from the Common 
 Core grade 5 material (5.NF(3)): 
\begin{quote}If 9 people
want to share a 50-pound sack of rice equally by weight, how many
pounds of rice should each person get?\end{quote}\medskip

\noindent Ratio approach: 
\begin{itemize} \item we have (50 lb rice) : (9 people)
\item we want (? lb rice) : (1 person) 
\item  1 = 9/9 so divide by 9 to get (\(\frac{50}{9}\) lb rice) : (\(\frac99=1\) person)
\end{itemize}
Ratios are thought of as two things that encode a relationship. They are not combined to give a single thing because it does not make sense to combine rice and people. Calculations in antiquity used methods that gave new pairs with the same relationship; nice illustrations are given by Schwartz \cite{palmiers} in a description of a fourteenth-century French manuscript. The division step in the example is a modern fraction-based shortcut applied to the numbers, and the pair-of-things format is retained to manage physical significance.\medskip

\noindent Rate approach: 
\begin{itemize} \item  (50 lb rice) = (\emph{rate})\(\times\)(9 persons)
\item divide to get  \emph{rate} =\(\frac{50}{9}\times \frac{\text{lb rice}}{\text{person}}\)
\item so for one person we get \((1\times\frac{50}{9})\,(\text{person}\times \frac{\text{lb rice}}{\text{persons}})\) = \(\frac{50}{9}\) lb rice. 
\end{itemize}
This reflects the seventeenth-century separation of physical magnitudes into numbers, and units or dimensions. Physical significance is handled by dimensional analysis that proceeds in parallel to the numerical manipulations. \medskip 

\noindent Discovery-learning approach:
\begin{itemize}\item Parse as ``blah blah 9 blah,  `share', blah 50 blah blah.''
\item answer: `share' \(\longrightarrow\) divide, so  \(\frac{50}{9}\)
\end{itemize}
The  material surrounding the 5.NF(3) example shows that students are supposed to infer how to work such things from examples, and are not given either explicit notation or logical clarity. But what they see is that almost all of the words in the problem statement are inert. Generally, problems are contrived, repetitive, and mathematically trivial: the rice problem is almost the same as examples given in 4NF(4c), 5NF(7c) and 6NS(1) and this is the standards document, not a text where repetition might be expected!  In practice the numerical data and operation required to get the answer can usually be extracted without reading the text, and this is the fastest and most reliable approach because it reduces cognitive interference between calculation and physical context; see the discussion in \S\ref{sssect:cogOverhead} below. 

\subsubsection{History} The ancient Greeks classified physical quantities as either magnitudes (continuous measures) or multitudes (discrete counts). They abstracted natural and rational numbers  from multitudes, and as a result these  numbers could be studied purely mathematically. The Greeks were not able to abstract real numbers from magnitudes, so continuous measures continued to have physical significance attached to them for another 2000 years. In effect every type of measurement had its own number system, and though they followed the same rules there were limits on how they could be combined because they had different `meanings'.  The number line, where rationals and irrationals live together in harmony, did not make sense because irrationals were not quite `numbers'. One consequence is that for 2000 years examples and exercises were \emph{necessarily} word problems and worked in context, because real numbers only made sense in physical contexts. Another consequence was the use of ratios as pairs of things,  thought of as an exemplar of a relationship rather than an object in its own right. The intended significance  of the relationship was inferred from context and not made explicit or precise. Not having to be explicit makes  physical applications seem easy, but this  is an illusion. It postpones rather than avoids the need for precision, and the lack of precision causes trouble in understanding numbers as well as in later physical applications. 

In the 1600s science professionals finally figured out how to decompose magnitudes into real numbers and units. 
By the mid to late 1800s the number/unit separation of magnitudes had evolved into `modeling', which separates entire physical problems into purely mathematical formulations, and physical significance. Analyzing the 
 mathematical component purely mathematically reduces cognitive interference and clarifies the role of mathematical structure. Analyzing the physical part qualitatively (for instance with `dimensional analysis') clarifies physical structure and relationships without distraction by complicated mathematics.  This separation has made  old problems easy and  more-complex problems accessible, and it has been standard practice in professional work for almost two centuries.  
 
\subsubsection{Cognitive overhead}\label{sssect:cogOverhead} The effectiveness of modeling has been clarified by developments in the cognitive sciences. We now know that humans have very limited working memory and very limited facilities for  logical analysis. We handle complex problems by breaking them into pieces, and in each piece focusing on essential features. Anything non-essential amounts to cognitive overhead that reduces our capacity to deal with the real problem. 

Mathematical analysis, and thinking about physical significance, are rather different activities that engage different neural regions and  strongly interfere with each other.
In particular the cognitive overhead in traditional word problems makes it very difficult for students to work them unless the mathematical core is trivial, or unless they learn to work them without engaging the physical context. Until grade 7 most problems have two or three numbers, and these are to be combined using addition, multiplication, etc. If taken at face value, most of the effort goes into deciding from the context which operations are to be used. The actual mathematics is clearly secondary. In practice these problems are so routine and repetitive that  most students learn to infer the operations without actually engaging the context. Some educational programs explicitly teach this with a `keyword' approach (``if you see `and' then  add, \dots''). In other words, the strategies students use to deal with this cognitive overhead often defeat the educational justification for word problems. But, unlike professional modeling, these strategies only work for mathematically-trivial problems and they do not open any doors. 

In particular, the physical-context format interferes with  understanding fractions. Division is sometimes not allowed, and it is not clear whether this is a physical-significance problem or a number problem. In very early grades the issue is 
dodged by evoking the  `parts-of-a-whole' loophole. Later it is 
partially addressed using the ratio and ``proportional reasoning'' point of view that keep different magnitudes separated. Vague  procedures and problem statements help hide the problem. All this reenforces the magnitude point of view. In later grades students are taught to get rate equations with fraction methods, and given examples to show how units are handled. But by then magnitudes are deeply embedded, and rate methods often come across as  arcane rules for the manipulation of magnitudes rather than triggering a change to the  the number-with-unit perspective. The lack of significant pure-number tasks, and the continued identification of word problems (where magnitudes seem to be needed) as the most meaningful and significant tasks, all contribute to the problem. An effective fraction perspective may never  emerge from the clutter. 

\subsubsection{Summary} Elementary education follows the ancient physical-context approach.  \emph{Experts} found this confusing, and abandoned it as counterproductive centuries ago. It is no wonder that students also find it confusing, and that many never recover from it.
\subsection{Dysfunctionality}\label{ssect:dysfunctional} 
Conceptual dysfunctionality due to the use of ancient cultural artifacts is essentally universal in elementary education. Consequences for skill development vary widely in different approaches.  The currently-dominant ``Reform'' movement has de-emphasized skills in favor of qualitative and subjective goals, aided by the use of calculators.  ``Traditional'' now means essentially ``not Reform'', and is  distinguished from it primarily by continued emphasis on functional skills. 

The Common Core Standards is essentially a Reform document. For instance, the eight ``Standards For Mathematical Practice" (immediately following the introduction)  set out ``varieties of expertise that
mathematics educators at all levels should seek to develop in their students''. Verbose abstract descriptions of this sort are  a characteristic feature of the Reform approach.  These descriptions sound more business-like than usual, probably because other groups participated in the development process, but they are vague, subjective, and deliberately designed to admit dramatically different interpretations. The sixth practice standard, for example, is ``Attend to Precision'', which to Traditional or professional readers sounds like  ``get the right answers''.  But the discussion goes in rather different directions and supports interpretations in which students could get high marks for ``precision'' \emph{without} getting right answers. 

To see how this might play out in practice we turn to another document.
\subsubsection{Multiplying fractions}
A group at the University of Pittsburgh designed a study to compare the effectiveness of the Reform and Traditional approaches in teaching the Common Core material. 
These alternatives were relabeled ``Dialogic'' and ``Direct (fully guided)'' to avoid some of the emotional baggage and, after convening several panels of experts (including the author) and soliciting feedback, the group developed working definition of the two approaches. The following comes from the ``Direct'' version \cite{pitt}:
\begin{quote}
\dots with respect to multiplying two fractions (e.g., \(2/3 \times4/5\)), fifth grade students would, in succession, need to:
\begin{enumerate}\item Know what a fraction is as a number (i.e., understand and identify the correct placement of fractions---whether posed in symbolic or non-symbolic form---on a number line);
\item Understand and use equivalence of multiple (symbolic) representations (e.g., knowing that \(2/3 = 2 \times1/3 = 1/3 + 1/3\)) on the number line---knowing, in particular, that all fractions can be decomposed into unit fractions (the ``Rosetta stone'' of fractions);
\item Employ commutativity (\(2\times 4 \times 1/3 \times 1/5 = 2/3 \times 4/5\));
\item Translate unit fractions (and operations over them) to the number line (\(1/3 \times 1/5\) on the number line is \(1/15\)); and
\item Use this to calculate a solution to the original problem: \(2/3 \times 4/5 = (2\times 4) \times (1/3 \times 1/5) = 2 \times 4\times1/15 = 8 \times 1/15 = 8/15\).
\end{enumerate}
\end{quote}
The subtext is that students should \emph{not} use the shortcut ``multiply tops and bottoms''. 
\subsubsection{Incomprehension} This snippet is most valuable for what it reveals about the Reform movement. A Traditionalist might say that skill development is a delicate business with many modes of failure. Fully guided instruction is needed to steer students away from these modes of failure and, when they occur, to detect and quickly correct them before they become so deeply embedded that they are long-term disabilities. In other words, instruction is driven by the desire to develop functional skills. For fraction multiplication, ``multiply tops and bottoms'' is the functional skill. The conceptual overlay in the snippet might justify why fractions are a good thing in general, or why they are appropriate in a particular word problem, but it is \emph{not} part of the skill. Why does this snippet so badly misrepresent Traditional goals? 

A cynical interpretation is that the authors are grafting Reform dysfunctionality onto the working definition of Traditional to sabotage it, so it will not look better in the study. I believe, however, that these are earnest and honorable people who simply cannot comprehend a skill orientation. Their preferred ``Dialogic'' approach allows students to make their own mistakes. They see that Traditional approaches use ``Direct (fully guided)'' instruction, but rather than understanding it as a way to \emph{avoid} mistakes they understand it as a way to get everyone to make the \emph{same} mistake. 
\subsubsection{Dysfunction} The procedure described in the passage puts heavy cognitive overhead on what should be a simple calculation, and clearly inhibits skill development. It is illuminating to see how this is accomodated by the Common Core.

The first of the Standards For Mathematical Practice, for instance, is ``Make sense of problems and persevere in solving them.'' Traditionalists might wonder how single-step word problems require perseverence. This is explained by the fraction-multiplication example: ``persevere'' refers to the elaborate mental gymnastics necessary to ``make sense of'', not the routine manipulation needed to ``solve''.  This shift of emphasis dilutes the importance of skills, and entitles teachers to say ``you can't work the problems but I can give you lots of partial credit because you `understand' the ideas.'' It is unfortunate that some of these ideas are mathematically incorrect and all of them have been purged from modern practice as counterproductive. 

What about cognitive overhead?  The text of step (5) shows that the writers are thinking of single-step word problems. There \emph{is} no further work that the overhead would interefere with. This reveals an interdependence in the Reform approach: cognitive overhead in word problems restricts the mathematical core to trivialities, often a single arithmetic operation. This enables de-emphasis of skills because even dysfunctional skills have a reasonable success rate with mathematically-trivial one-step problems. Especially when calculators are used. Dysfunctional skills then lock in the committment to simple word problems, because students have trouble doing more. A strong word-problem orientation is therefore essential for the Reform approach to appear to work at all. This clarifies the word-problem orientation of the Common Core Standards. 

\subsubsection{Summary} The modern expert-user version of fractions has evolved to maximize functionality, driven by the needs of science. The    fraction material in the Common Core Standards is heavily influence by ancient cultural artifacts, some of them predating mathematics-based science by two millennia. Some of it is mathematically incorrect, and even the correct parts are not particularly functional. This material is common to both traditional and Reform approaches to elementary education. Traditional teachers usually try to develop some functional skills in spite of the material, while Reform educators emphasize ancient ``meanings'' over skills. The Common Core Standards support the de-emphasis of skills, for instance by expanded focus on mathematically-trivial word problems. 

\section{Mathematics of fractions} \label{sect:carefulfrax} The goal is a full-precision description of how mathematicians use fractions, and of the underlying mathematical structure.  Subtleties have to be avoided when such material is taught to children, but it seems to me this means educators should have complete clarity about these subtleties. Without this---as we have seen---children will get stuck in old dead ends and their learning will not be upward-compatible. 

Fractions and negative numbers are considered together because the underlying mathematics,  motivations, and patterns of  expert use are almost identical. Incidently, this is why mathematicians find fractions exactly as easy as negatives, and why schoolchildren could too.
 
\subsection{The setting} The setting is a commutative semiring. This means multiplication and addition are defined and satisfy the familiar rules of arithmetic, but negatives and quotients are not always defined.  Examples encountered in K-12 include natural numbers, non-negative real-valued functions, and polynomials with nonnegative coefficients.  
When the focus is specifically on fractions we will work in a commutative ring (not semi-). The difference is that elements are assumed to have additive inverses. 
\subsection{Functionality}\label{sssect:mainjob} The primary job of negatives and fractions is easy manipulation of equations. The simplest example is the description of solutions of first-order  equations as a formula in the coefficients. To be explicit, suppose \(a\times\square +b=c\). Unadd  to get \(a\times \square =c-b\), then unmultiply to get the formula \(\square =\frac{c-b}{a}\). \medskip

The terms ``unadd'' and ``unmultiply'' are used to emphasize that (unlike subtraction and division) they are \emph{operations}. Invoking these operations may cause a change of rings to make them defined, but for the most part this need not concern students. The machinery that makes this work is described in \S\ref{ssect:fraxrings} and educational implications in \S\ref{ssect:edimps}; here I explain what  is important about it, and why. 
\subsubsection{Cognitive simplicity}
 The operational advantage of the unadd/unmultiply view is cognitive simplicity. We start with \(a\times\square +b=c\) and want to get \(\square\) all by itself on one side of an equation. The \(b\) is in the way so we want to move it to the other side. This is an organizational decision, and separating organization from numerical or algebraic activity  reduces cognitive interference and overhead. It is therefore a considerable advantage to move the \(b\) with a formal operation: ``unadd \(b\)\dots''.  
 
After  \(b\) has been moved we see the coefficient \(a\) as being in the way, and make a organizational decision to move \emph{it}. Again this is accomplished with a formal operation, ``unmultiply to move \(a\) to the other side''. These operations are the same for numbers, symbols, polynomials, functions, etc., so do not engage calculational activity. Only after the organizational phase is complete do we shift gears to process the resulting expression. 

Note that  the procedure has to be \emph{described} as a formal operation to get the full benefit. The usual ``subtract \(b\) from each side'' or ``divide each side by \(a\)'' requires cancellation of \(b-b\) or \(a/a\) on the left as an additional step. This wastes the information that the operation is \emph{designed} to cancel. 
The usual approach also mixes organizational activity with arithmetic so  is slower and increases the error rate. Finally, because it requires a ``shift of gears'', it can be a serious distraction. When this procedure occurs as a step in a larger problem I often see students' larger-scale train of thought disrupted by such distractions. 

\subsubsection{Efficiency driven by science} Formal operations give cognitive simplicity, but for 2000 years people got along pretty well without it and in fact without negative numbers at all, and with ratios instead of fractions. Why is this a big deal now?

The job of mathematics itself changed radically in the 1600s with the development of science based on serious mathematics. Newton's laws were just ink on paper until mathematics brought them to life, but then they became so powerful that they changed the world. This put a lot of new pressure on mathematics. It took a while to develop the mathematics to implement Newton's laws and later developments follow the same pattern.  When Maxwell developed a unified description of electricity and magnetism in the 1860s for instance, it was almost useless because it was beyond the reach of the mathematics of the time. It  became powerful when mathematics caught up. The same pattern occurred with relativity, quantum theory, and many other developments. 

The point is that by the late 1600s ``advanced'' mathematics had gone from being a relaxed intellectual exercise to the engine of science. Ambitious goals and urgent demands for progress came from outside the community. Things like solving first-order equations went from occasional exercises to ubiquitous components of larger problems. Large-scale use drove evolution of  versions optimized for human use, much as the demands of commerce drove the optimization of arithmetic algorithms long before. The cognitive simplification described above was a major part of this optimization.

I mention two educational consequences. The first is that by 1700  essentially all K-12 math, and fractions in particular, had become tools rather than destinations. Professional users needed mastery  of symbolic versions optimized for accuracy and transparency, and had no use for the old ``meanings'' or visualization of analogs or special cases. Elementary education does not reflect this change, and consequently does not prepare students for the technical careers of even the nineteenth century. 

The second consequence is more subtle. 
Methodological evolution driven by heavy use is often invisible to the users. They may know why it works, but are not consciously aware of features that make it work \emph{well}. Such methods  cannot be taught explicitly because teachers do not realize there is something to teach. Instead, students get started by imitating teachers, then develop full skill through heavy use and a personal micro  version of methodological evolution. It follows that watching an experienced role model work through material teaches much more than generally realized. A sample consequence is that computer-based courses must be supplemented \emph{at least} by videos of experienced users working problems in order to be fully successful. This also clarifies that ``experienced'' means ``enough heavy use of the methods to develop automatic use of optimized versions''. In other words, effective teaching of subconscious knowledge requires  teachers who are over-qualified in the subject matter. Expertise in \emph{teaching} by itself has no benefit.
\subsubsection{Constraints}\label{sssect:constraints}  Formal operations can be powerful, but there are limitations. Successful human use depends having few problematic cases (eg.~dividing by zero) and knowing exactly what goes wrong. Non-commutative rings (eg.~square matrices) are sufficiently more complicated that we speak respectfully of ``multiplying by inverses'' rather than either division or fractions. The infinitesimal fraction approach to calculus is too tricky for ambitious use (see \S\ref{sssect:oldcalculus}). And sometimes negatives or fractions collapse the number system in unacceptable ways (see \S\ref{sssect:tropical}). 
The point is that effective mathematical methodology is an evolved balance between human limitations and mathematical structure. History, philosophy, and wishful thinking are poor guides.
\subsection{Universal objects}\label{ssect:fraxrings} We turn to the technical details. Examples are described in the next section, and implications for education in the section after that.

  Suppose we have an associative and commutative operation \(\mathsf{op}(a,b)\) on a set \(R\), and a collection \(D\) of elements for which we want to have an
  ``un-operation'' \(\mathsf{unop}(d,a)\) satisfying  \(\mathsf{op}(d,\mathsf{unop}(d,a))=a\), if \(d\in D\).  
  The standard way to do this is to consider operation-preserving morphisms to other sets-with-operation, \(R\to \hat R\) so that
  \begin{enumerate}\item \(\hat R\) is a  has a neutral element \(e\) for the operation, and
  \item for every \(d\in D\) there is a unique \(q\in \hat R\) so that \(\mathsf{op}(\hat d,q)=e\).
  \end{enumerate}
 In any such \(\hat R\) we can define \(\mathsf{unop}(d,a)=\mathsf{op}(q,a)\). 
 
 The next step is to consider morphism which are universal with respect to this property. Explicitly, \(R\to U\) is \emph{universal} if given \emph{any} morphism  \(R\to S\) such that the desired un-operations exist in \(S\), then it factors through a unique morphism \(U\to S\). 
 There is a simple and very general argument (called ``general nonsense'' by practitioners) that implies there are  morphisms that are universal for this property, and that these are all equivalent. 
 
\subsubsection{The key lemma}  There is a standard explicit description of a universal object: start with the set of ordered pairs \((a,b)\) with \(a\) obtained by a finite iteration of operations using elements of \(D\), and divide by an appropriate equivalence relation.  Details (omitted here) obscure the cosmic inevitability of universal objects but make it easy to show a key technical result:
\begin{description}\item[Lemma]  If all elements of \(D\) satisfy cancellation (\(\mathsf{op}(d,x)=\mathsf{op}(d,y)\Longrightarrow x=y\))  then  universal morphisms are injective. Conversely, if \(D\) contains an element that does not satisfy cancellation then  universal  morphisms are not injective.
\end{description}
Ask a mathematician about fractions and you will probably get some version of this explicit description. As observed above this is the technical key, but does not give insight into motivation or useage.

\subsection{Examples of negatives and fractions}\label{ssect:constructfrax}  
These illustrate the uses and limits of fraction-like constructions with commutative operations. The fraction approach to calculus described in \S\ref{sssect:oldcalculus} is particularly illuminating because it seemed to work pretty well for a while, but eventually---to the dismay of naive users---had to be replaced with limits. 

\subsubsection{Natural numbers}\label{sssect:unadd} Let \(N\) denote the natural numbers \(1,2,\dots\). Addition and multiplication both define  commutative and associative operations, and every element satisfies cancellation with respect to both operations.  According to the `key lemma'  \(N\) injects into universal monoids where the operations can be inverted. Inverting multiplication gives positive rationals, and inverting addition gives the integers. These are the standard descriptions in modern algebra courses. 

These two constructions are not quite symmetric. \(N\) contains a neutral element for multiplication, and it follows that elements obtained by inverting multiplication still satisfy additive cancellation. Inverting addition then gives  rationals numbers. \(N\) does not contain a neutral element for addition (zero). Inverting addition introduces one, but it does not satisfy multiplicative cancellation. Wholesale inversion of multiplication as a second step therefore gives the zero ring. To avoid collapse we invert only non-zero elements.

Historical comments: The ancient Greeks  abstracted the natural numbers and used un-multiplication to expand them to the positive rationals. Expansion to include negative numbers, essentially by un-addition,  seems to have been accomplished in India by  630 CE \cite{Brah}, but they did not know what to do with the failure of 0 to satisfy multiplicative cancellation (ie.~the consequences of zero denominators). They were still fudging this 600 years later. The Arabic world got base-ten representations of numbers from India but seem to have found negatives too problematic to adopt at the same time. One reason is that ratios were used as the standard division-like operation in the west until the 1600s, and negative numbers cause much more trouble with ratios than with fractions; see \S\ref{ssect:carefulratios}. 

The use of unaddition to modify  monoids to have additive inverses was formalized around 1950 (?), and is called the ``Grothendieck construction''.   Fancy name for a particularly easy version of fractions. The cancellation condition is not required so typically there is some collapse. Usually this is thought of as a good thing: what is lost is subtle ``unstable'' information that interferes with a coherent global view. Sometimes it is a bad thing, as we see next. 
\subsubsection{Tropical arithmetic}\label{sssect:tropical} The \emph{tropical semiring}\footnote{See the Wikipedia entry.} is the real numbers and \(+\infty\) with the operations 
\begin{itemize}\item \(r\oplus s =\min(r,s)\), and
\item \(r \otimes s = r+s\).
\end{itemize}
The notations \(\oplus, \otimes\) reflect the fact that both operations are commutative and associative, and that \(\otimes\) distributes over \(\oplus\). Looks like a ring, but rings are usually required to have inverses for  ``addition''  and this doesn't.  We can't do anything about it either: ``additive'' cancellation fails for every element except \(+\infty\). Specifically, if \(r,s\geq a\) then \(r\oplus a=a=s\oplus a\). Requiring un-\(\oplus a\) to be defined therefore puts us in a semiring where the whole interval \([a,+\infty]\) collapses to a point. 

\subsubsection{Polynomial fractions} If \(P, Q\) are real-coefficient polynomials then \(P/Q\) is usually called a `rational function'. This is a misnomer because we use fraction methods to work with them, and we use the fraction meaning for equality of two of them. To illustrate this last point consider
\[\frac{a^2-x^2}{a-x}\stackrel{?}{=}a+x\]
This is true for fractions, but false for functions because one side is defined when \(x=a\) and the other is not. 

 \subsubsection{Formal calculus}\label{sssect:oldcalculus}  The differential calculus of Newton and Leibniz  is essentially  fraction algebra. Why it worked for a while and then failed illustrates the limits of the fraction construction. 

Suppose \(f\) is a power series in \(x\), without any assumption of convergence. We may think of \(f\) as having real coefficients, but the discussion works for any commutative ring.  

If \(x\in R\) then \(f(x+t)\) can be written as a  series in \(t\). The constant term is \(f(x)\) so  \(f(x+t)-f(x)\) is divisible by \(t\). Define the \emph{derivative} \(Df(x)\) to be the constant term in the quotient  \(\frac{f(x+t)-f(x)}{t}\). 

Note that writing \(f(x+t)\)  as a  series in \(t\) implicitly  extends \(f\) to a series defined on the polynomial ring \(R[t]\). In this context the constant term of a series is the image under the morphism \(R[t]\to R\) that takes \(t\) to 0. This shows it is absolutely essential to be sure that \(f(x+t)-f(x)\) is divisible by \(t\) in  \(R[t]\). If it is not, then writing \(\frac{f(x+t)-f(x)}{t}\) takes us to something like \(R[t, \frac1t]\), and then setting \(t=0\) takes us to the zero ring rather than \(R\). All without warning. 

This can be streamlined to avoid the ``set \(t=0\)'' step by using the quotient 
ring \(R[\dt]/<\dt^2=0>\). Elements are of the form \(a+b\,\dt\), so we could also think of them as pairs \((a,b)\) with multiplication \((a,b)*(r,s)=(ar,as+br)\). Notice that \(dt\) does \emph{not} satisfy cancellation so forming fractions with this in the denominator is problematic.  Elements with invertible real part are invertible however, and the inverse is: 
\[\frac1{a+b\,\dt}=\frac1{a+b\,\dt}\times\frac{a-b\,\dt}{a-b\,\dt}=\frac{a-b\,\dt}{a^2}=\frac1a-\frac{b}{a^2}\,\dt.\]

Again, a power series on \(R\) extends  to  \(R[\dt]\) simply  by plugging in elements of the extended ring. The  extended function works out to be
\[f(x+y\,\dt)=f(x)+D_{x}f\cdot y\,\dt\] with  \(D_xf\) as above. 
Note that applying this to \(f(x)=x^{-1}\) gives the formula for inverses displayed above.\footnote{Like some other fraction formulas this one looks a bit odd. The derivation shows that it is formally correct but some students and some educators may see this as a sign that the whole endeavor is artificial. Connecting the formula to the derivative of the inversion function is a mathematical-understanding maneuver: the formula looks odd in isolation but reasonable and natural in a richer context.  }

Next, set \(y=1\) in this expression to get \(f(x+\dt)-f(x)=D_{x}f\cdot \dt\). Dividing by \(\dt\) then gives the seventeenth-century description of the derivative as a fraction: \[D_xf=\frac{f(x+\dt)-f(x)}{\dt}.\]
This maneuver makes the seventeenth-century formula sensible and correct for power series but, as noted above, it is problematic because \(\dt\) does not satisfy cancellation. If we ever have a function with an extension to \(R[\dt]\)  such that \(f(x+\dt)-f(x)\) is  \emph{not} divisible by \(dt\), then the fraction  puts us in the 0 ring. 

 Beginners had trouble with this approach because it  can fail disastrously if \(\frac{b}{\dt}\) is manipulated too much like a genuine fraction. Experts learned where the edges of validity were and stayed away from them, but could not avoid trouble when they began to work with non-analytic functions. Eventually this approach was replaced by the use of limits because limits are more general and more robust. 
 
It is useful to understand \emph{why} limits are more robust. Limit considerations start with a presumption of \emph{non}-existence, and any use should include a justification of  existence. In principle, mathematical arguments are \emph{supposed} to be justified. It is a lot of trouble, and occasionally there is a perfectly good limit that we can't use because we can't justify it, but in return we get full assurance that what we do makes solid sense. Fractions start with a presumption of \emph{existence} and are much easier to use because we don't have to justify anything. This is a violation of basic mathematical procedure because arguments can become nonsense without warning if we overstep the bounds of validity. When the bounds are very simple (``don't divide by zero'') humans can do this and stay out of trouble, but limits can fail in far too many subtle ways for us to be able to use this approach safely. 
\subsubsection{Non-standard fractions} In the 1960s Abraham Robinson gave another justification of seventeenth-century calculus using what amounts to a variation on the fraction construction. Suppose \(d\) is a ring element that does not satisfy cancellation (think \(d=\dt\) as above). Writing \(\frac1d\) puts us in another ring where it does satisfy cancellation, so an easy argument shows something has to be collapsed. Robinson discovered that sometimes you can avoid collapse by making the easy argument \emph{illegal}! Technically, \(d\) and \(\frac1d\) are in a first-order logical extension of the real numbers. If they are used in  appropriate first-order logical expressions then the expressions have sensible output without forcing the ring to collapse. However accidental or careless use of an overly powerful (second-order) expression may cause the number system to collapse without warning. This is dangerous and hard to manage. There have been attempts to package this approach to hide the logical complexity \cite{hyperreals}, but again the limit approach is more flexible and robust.

\subsection{Educational implications} \label{ssect:edimps} We can only scratch the surface of this topic.  The full-precision description of fractions is not something to teach: no-one can think it would be a good idea to teach third-graders to say ``universal object''. Instead it should  guide development of context and subliminal influences in teaching. I give a few illustrations.

\subsubsection{Automatic contexts} Fractions determine their own context in the sense that writing  \(\frac{b}{a}\) puts us in a ring where the quotient (of images) exists.  Writing \(\frac32\) automatically puts us \emph{at least} in \(Z[\frac12]\), the closest approximation to the integers in which the quotient \(3\div 2\) is defined. Our notation does not record the change of rings, but this is harmless as long as the change is injective. The net result is that for elementary use it is unnecessary to say anything about \emph{what} or \emph{where} a fraction is, as long as we avoid dividing by zero. 

The corresponding drawback is that \emph{attempting} to explain the `what' and `where' of fractions is likely to turn a straightforward procedural subject into a confusing mystery. The precise details show that writing a fraction specifies an object in the image of \emph{every} appropriate ring morphism, not just one. One might think of this as being like a picture on the internet: it doesn't appear \emph{every}where; one must go to a device with an internet connection, but it appears in every such context. This is a perfectly functional way to think about internet pictures. A full-precision description as functions, from  internet-connected devices and URLs to images, would not help anybody. Falsehoods analogous to the ancient descriptions of fractions are not helpful either because they interfere with development of genuinely functional understanding. 

Similar considerations apply to special collections of fractions. The full-precision view of rational functions, or rational numbers, is that the essence is in \emph{specifications of structure}. Specific implementations are called \emph{models} for the structure, and there are lots of them. One of the consequences of Russell's paradox is that none of them can claim to be \emph{the} rational functions or \emph{the} rational numbers. The ordered-pair model is useful but does not have a privileged status.  The models are all equivalent so there is nothing wrong with talking about ``the rational numbers'' as long as one remembers that the word ``the'' is a linguistic artifact not to be taken literally. 
 
\subsubsection{Zero denominator}\label{sssect:zerodivisors} Writing \(\frac20\) automatically puts you in the zero ring (ie. \(R[\frac10]=0\)).  After this, like it or not, everything takes place in the zero ring, so \emph{everything} is zero. We avoid putting 0 in denominators to avoid causing the number system to collapse, not because it is undefined. The ancient admonition ``don't divide by zero'' works fine as a procedural rule and does not need to be changed, but the traditional explanation ``not allowed'' does need to be changed. 

Note that `if it isn't broken, don't fix it' applies to this. Students who ask about division by zero may be receptive to such an explanation. Students who \emph{don't} ask will probably not understand why it is an issue and are more likely to be mystified than enlightened. 

I give some historical perspective. If you assume \(\frac{2}{0}\) is defined then you can easily show that \(1=0\). This seems obviously false. The old response follows the patterns of science: any physical theory has limits of validity and if you go past the edge you fall off and get hurt. You learn where the edges are and stay away from them.  Dividing by zero seems to be beyond the edge of the world of fractions. This was a fair description until the late nineteenth century when this sort of thing became intolerable. Mathematicians work very close to edges like this and have to know exactly where they are and exactly what goes wrong.  There were determined efforts in the nineteenth and early twentieth centuries to  resolve all such paradoxes and, judging by a century of ultra-extreme testing since then, the efforts were successful. How was the  \(\frac{2}{0}\Longrightarrow (1=0)\) problem resolved? Careful precision reveals that the meaning of ``='' depends on context, and \(\frac{2}{0}\) puts us in a context where  \(1=0\) is indeed true. The mistake is to assume that ``='' means what you want it to mean. The implicit presumption that the original ring injects into the fraction ring is wrong in this case. This is a defect in understanding, not a defect in mathematics. Zero denominators are still something we avoid, but now we know exactly why. 

\subsection{Mathematics of ratios}\label{ssect:carefulratios} Ratios can be developed in general commutative rings, but the main lesson is that one has to work harder to  get a whole lot less. We skip this and go directly to ratios of real numbers. 

\subsubsection{Parameterizing the plane} Think of pairs \((a,b)\) of real numbers as points in the plane, then every pair except \((0,0)\) determines a ratio, and every pair with nonzero first coordinate determines a fraction. The clearest picture comes from the inverse process: using ratios and fractions to parameterize the plane.

For both fractions and ratios, the points equivalent to \((a,b)\) lie on the straight line through  \((a,b)\)  and the origin. We can therefore describe a point in the plane as the line (= all equivalent fractions or ratios) , together with a point in the line. This parameterization is injective except that, since the lines pass through the origin, the zero point in every line goes to the zero in the plane. The ratio parameterization gives the whole plane, the fraction parameterization gives everything but the vertical axis. 

In mathematics the construct ``line, together with a point in the line'' is called the canonical line bundle over the projective space. The mathematical problem with ratios is that the canonical line bundle is a M\"obius band: it has a twist. 
I mentioned in \S\ref{sssect:unadd} that negative numbers  developed early in India but did not catch on in the Arabic and western world  because they interact poorly with ratios\footnote{This was still a problem in 1600: Descarte found negative numbers so problematic (for ratio reasons) that he referred to them as ``false numbers''. Complex numbers seemed relatively harmless so were merely ``imaginary''.}. The bad interactions are tedious to illustrate explicitly, and in particular cases one might imagine that they could be resolved by being more clever. But they have their origins in the twist in the M\"obius band, so the ratio approach is doomed to be twisted no matter how clever we are. 

In these terms we see fractions as  avoiding the twist by cutting the M\"obius band along the preimage of  the vertical axis in the plane.  Another way to avoid the twist is to restrict the equivalence relation  \((a,b)=d\times(x,y)\) to \emph{non-negative} \(d\). This cuts the band along the central circle and gives an oriented band that double covers the original. The corresponding parameterization of the plane is polar coordinates.  

\subsubsection{Summary} Ratios are, by nature, twisted. Orienting them to  avoid the M\"obius twist gives angles in polar coordinates, not fractions. 
\section{Diagnosis} We have followed a single thread (fractions),  mostly through a single document. Following other threads gives other insights, particularly about some of the cognitive issues, but so far they all lead to the same general picture. 

Human methodologies evolve in response to selective pressures, and most rapidly in response to the greatest pressures. For about four centuries mathematical methodology has been driven by the needs of science and, especially in the last century, increasing accessibility of ambitious goals within mathematics itself \cite{rev}.  Modern practice is immensely more powerful as a result. 

Elementary education was largely insulated from the pressures driving the profession. Instead, for at least the last century, the main pressure on educational methodology comes from the need to sell it to administrators, legislators, and the general public. The current Reform movement represents a breakthrough in this direction, as powerful in its own way as any technical innovation in professional practice.  
 They have made a \emph{reduction} in skill levels sound exciting, and done it so well that they have become the dominant movement in just a few decades. The introduction to the Common Core document, and the  ``Standards For Mathematical Practice" that follow, offer visions and abstract goals that are much more compelling than any Traditional or Modern account. They have the practical effect of reducing skill expectations almost to zero, but this gives the Movement a further tactical advantage: modern methods are distinguished by their efficiency and power, but in the Common Core there is not much for them to do. A hammer does not look good when it is used to squash ants. 
 
Perhaps, someday, children will experience the clarity and power of modern mathematics. We are currently moving in the opposite direction, however, and it is hard to see how this could change.  

\enddocument